\documentstyle[leqno]{article}

\newcommand{\ol}{\overline}
\newcommand{\bqq}{\begin{equation} \label}
\newcommand{\eeq}{\end{equation}}
\begin{document}

\renewcommand{\thefootnote}{\fnsymbol{footnote}}
\begin{titlepage}

July 1995\hfill{KSU/GRG-95-4.7}
\vskip3cm
\begin{center}
{\bf QUANTIZATION OF K\"AHLER MANIFOLDS ADMITTING $H$-PROJECTIVE
MAPPINGS}\\
\vskip2cm

A.V.Aminova \footnote{E-mail: aminova@phys.ksu.ras.ru},
D.A.Kalinin \footnote{E-mail: kalinin@phys.ksu.ras.ru}\\
Dept. of General Relativity and Gravity, \\
Kazan State University, Lenin Str., 18, \\
KAZAN, 420008, Russia
\vskip1cm

\end{center}
\begin{abstract}
We discuss the  quantization  of  mechanical  systems  for  which  the
Hamiltonian  vector  fields  of  observables  form  the deformation of
$n$-dimensional oscilator algebra.  Because of this fact these systems
can  be  considered as "deformations" of the harmonic oscillator.  The
set of above-mentioned mechanical systems are realized at the
classical
level  in  the  form  of  K\"ahler  manifolds  of constant holomorphic
curvature.  Such mechanical systems are quantized later with the  help
of   the   geometric
quantization  approach.  We  also  discuss  the
quantization of more general K\"ahler manifolds  (not  necessarily
of constant holomorphic curvature) admitting $H$-projective mappings.
\end{abstract}
\end{titlepage}
\renewcommand{\thefootnote}{\arabic{footnote}}
\setcounter{footnote}0
\setcounter{section}{0}

%%%%%%%%%%%%%%%%%%%%%%%%%%%%%%%%%%%%%%%%%%%%%%%%%%%%%%%%%%%%%%%%%%%%%%
\section{ Introduction}

Let us call the
$n$-dimensional oscillator algebra ${\bf osc}(n)$
Lie algebra with
$(2n+1)$ basic  generators  $T^{\alpha},T^{\ol\alpha}$  and $T$
which obey the following commutation relations
$$
[T^{\alpha},T^{\ol\beta}]=[T^\alpha, T^\beta]=
[T^{\ol\alpha},T^{\ol\beta}]=0,
$$
\begin{equation} \label{te1}
[T^{\alpha},T\hskip 1pt ]=-i\hskip 1pt T^{\alpha}, \qquad
[T^{\ol\alpha},T\hskip 1pt ]=i\hskip 1pt T^{\ol\alpha},
\qquad \alpha=1,...,n.
\end{equation}

Let ${\cal  A}_m$  be  the  manifold  of  Lie  algebra  structures  in
an $m$-dimensional vector space $V_m$. The curve $ \ell (t) $
in  ${\cal A}_m$ passing through the algebra
$\ell =\ell (0)$ is called the deformation  of  the  Lie  algebra
$\ell$.

It is   well   known   that   deformation  of  one-dimensional
oscillator algebra  play  an   important   role   in   various
considerations. They are
closely related   with   many   interesting    physical    and
mathematical objects,  for example, anti-de Sitter quantum mechanics
[{\bf 4,6}], symplectic geometry of one-dimensional complex disk
[{\bf 5,7,14}] and quantization of spinning particle [{\bf 12}].

This paper is devoted to quantization of the mechanical systems
for which the  Hamiltonian  vector  fields  of  observables  form  the
deformation of $n$-dimensional oscilator algebra. Because of this fact
these systems  can be considered as "deformations" of
the harmonic oscillator. The paper consists of
three sections.  In the first section we give short review of the
geometric
quantization procedure. In the second section is
constructed  the  set  of above-mentioned
mechanical systems which are  realized  at  the
classical level in the form of K\"ahler  symplectic  manifolds
of constant   holomorphic   curvature  [{\bf 9}].  Such  mechanical
systems
are quantized later with the help of the geometric  quantization
approach [{\bf 8,12}]. As  it is  known  the  K\"ahler manifold  of
constant holomorphic  curvature   is   $H$-projectively   flat,
i.e. it admits $H$-projective mappings on the flat space. In
the third section of the paper we  discuss  the  quantization  of
more general  K\"ahler  manifolds (not necessarily of constant
holomorphic curvature)  admitting  $H$-projective  holomorphic
mappings.
\\

%%%%%%%%%%%%%%%%%%%%%%%%%%%%%%%%%%%%%%%%%%%%%%%%%%%%%%%%%%%%%%%%%%

\section{ Preliminaries}

To start with, we  recall  some relevant
facts   about   geometric
quantization procedure [{\bf 8,12}]. Let $(M,\omega)$ be a symplectic
manifold. According   to   Dirac   {\it quantization} is the linear
map
$Q:f\rightarrow \hat f$ of Poisson (sub)algebra $C^ {\infty} (M)$ into
the set  of  operators in some (pre)Hilbert space
${\cal H}$ possessing the properties:

\begin{enumerate}

\item $\hat 1=1$;

\item $\widehat {\{ f,g\}}_h = \frac{i}{h}(\hat f \hat g - \hat g \hat
f)$;

\item $\hat {\ol f} = (\hat f)^{*}$;

\item for  some complete set of functions
$f_1 ,...,  f_n$ the operators ${\hat  f}_1 ,...,{\hat f}_n$ also
form the complete set,
\end{enumerate}
where bar  is  the  complex  conjugation, star
denotes  the conjugation of the operator and $h= 2\pi\hbar$
is the Planck constant.
The linear  map  ${\cal  P}:f  \rightarrow \check f$ possessing the
first
three properties   is   called   {\it prequantization}. For the case
$M=T^{\ast}M$, $\omega = d \alpha $ prequantization was constructed
by Koopman, Van Hove and Segal. It has the form
\bqq{te2}
 {\cal P}f=\check f = f-i \hbar V(f) - \alpha (V(f)),
\eeq
where vector  field
$V(f)$ is the Hamiltonian vector field of the function $f \in
C^{\infty} (M)$ defined by the condition
$$
V(f)\hskip 1pt \rfloor \hskip 1pt \omega = -df,
$$
where $\rfloor$ is the internal product.  In local  coordinates
$x^i$, $i=1,...,2n$ from here we have
\bqq{te3}
V(f) = \omega ^{ij}\partial_j f\partial_i  .
\eeq

Let ${\cal L}$ be the Hermitian line bundle with  connection  $D$
and $D$-invariant   Hermitian   structure   $<,>$.   Recall  that
$D$-invariance means that for each pair of sections $\lambda$ and
$\mu$ of ${\cal L}$ and each real vector field $X$ on $M$ holds
\bqq{te4}
X<\lambda ,\mu  >=<D_X  \lambda  ,\mu>+<\lambda,D_X  \mu>.
\eeq

Let $(x,U)$ is local coordinate system on $M$.  If $\mu_0$ is
nonvanishing section of ${\cal L}$ over $U$,  then we can identify the
space $\Gamma  ({\cal  L},  U)$  of section with $C^{\infty}(M)$ by
the
formula $C^{\infty}(U)\owns \varphi \leftrightarrow \varphi \mu_0 \in
\Gamma
({\cal L},  U)$.  The  operator $D_X$ in this case takes the next
form
\bqq{te5}
D_X \varphi = X\varphi -i\hbar^{-1}\alpha(X)\varphi,
\eeq
where 1-form $\alpha$ is given by the relation
\bqq{te6}
D_X \mu_0 = -i\hbar^{-1}\alpha(X)\mu_0 .
\end{equation}
By comparing (\ref{te2}) and (\ref{te6}) we have the prequantization
formula of
Souriau-Kostant
\bqq{te7}
\check f = f-i\hbar D_{V(f)}.
\end{equation}
The curvature form $\Omega$ of $D$ is defined by the identity
$$
\Omega(X,Y) =     \frac{1}{2\pi i}([D_X,D_Y]       -
D_{[X,Y]}),
$$
and locally we have
\bqq{te8}
\Omega = h^{-1} d\alpha.
\end{equation}

{\bf Theorem  1}  [{\bf 8}].  {\it The  Souriau-Kostant  for\-mu\-la}
(6) {\it  defines
the pre\-qu\-an\-ti\-za\-ti\-on
if and only if  the  curvature  form  $\Omega$
coincides with} $h^{-1}\omega$.
\\

The construction of Hilbert space ${\cal H}$ in
geometric  quantization  procedure
essentially involves  the  choice  of the {\it polarization} that is
the involutive   Lagrange  distribution   $F$
in  $TM\otimes_{\bf R}{\bf C}$.
The polarization $F$ is called {\it K\"ahler}  if
$F\cap\ol F=\emptyset$,   and   the  Hermitian  form  $b(X,Y)=
i\hskip 1pt \omega (X,\ol Y)$ is positively defined for $X\in F$.

If the polarization $F$ on $TM\otimes {\bf C}$ is chosen then
the Hilbert space ${\cal H}$ consists of the sections
$\lambda$ of ${\cal L}$  which  are covariantly constant along $F$
\bqq{te9}
D_X \lambda=0,\quad\lambda\in{\cal L},\quad X \in F.
\end{equation}

It is  said  that  the  function  $f  \in C^{\infty}(M)$ preserves the
polarization $F$
if the flow $V(f)$ of $f$ obeys the condition
\bqq{te10}
[V(f),X^{\alpha}]=a[f]^{\alpha}_{\beta}X^{\beta},
\end{equation}
where $a[f]^{\alpha}_{\beta}$ are  some smooth functions  on  $M$
and vector
fields $X^{\alpha}$, $\alpha = 1,...,n$ span  $F$.
Let us  consider  the  particular  case  when  the polarization $F$ is
spanned by the complex Hamiltonian vector fields. In this case the
functions preserving
polarization can be  quantized  with  the  help  of  the  next
formula [{\bf 12}]:
\bqq{te11}
Qf=-i\hbar D_{V(f)}+f- \frac{i\hbar}{2}a [f],
\end{equation}
where $a[f]=\Sigma^{n}_{\alpha=1}a[f]^{\alpha}_{\alpha}$.

In general situation $f$ it does not preserve the po\-la\-ri\-za\-ti\-
on
and for qu\-an\-ti\-zing $f$ one have to use the
{\it Blat\-tner-Kos\-tant-Stern\-berg (BKS) kernel} which    connects
representations for different po\-la\-ri\-za\-ti\-ons
[{\bf 8,12}].  It is easy
to see that ever in the case of an $n$-dimensional  oscillator
the flow of the Ha\-mil\-tonian
\bqq{te12}
H=\frac{1}{2}\sum\limits^{n}_{\alpha =1}((p^{\alpha})^{2}
+ (q^{\alpha})^{2})
\end{equation}
does not preserve  the  polarization  spanned  by the  Hamiltonian
vector fields  of  both
position $q^{\alpha}$  and momentum $p^{\alpha}$ variables.
Therefore we must use BKS kernel to quantize $H$.  However,  when
we introduce the complex coordinates
$$
z^{\alpha} =   \frac{1}{\sqrt{2}}(p^{\alpha} + iq^{\alpha}),
\qquad
z^{\ol\alpha} = \frac{1}{\sqrt{2}}(p^{\alpha} - iq^{\alpha})
$$
and K\"ahler   polarization   spanned  by  the  vector  fields
$V(z^{\alpha})$ we   can   see   that   $H=\Sigma   z^{\alpha}
z^{\ol\alpha}$ preserves the polarization
and one  can use (\ref{te11}) to quantize $H$.  As the result they
obtain
the differential operator $\hat H$
\bqq{te14}
\hat H= Q_{\cal FB}\hskip 2pt H  =
\hbar (z^{\alpha}\frac{\partial}{\partial z^\alpha}+
\frac{n}{2})
\eeq
on the space ${\cal O}(U)$ of holomorphic functions  on
$U\in {\bf C}^{n}$ which we denote $Q_{\cal FB}$  because  the
corresponding rep\-re\-sen\-ta\-tion     is    called    Fock-Barg\-
mann
rep\-re\-sen\-ta\-tion.
In  the considered  case the rep\-re\-sen\-ta\-tion space ${\cal H}$
consists of the sections of ${\cal L}$  which  have  the  form
$\psi (z) \mu_{0}$,  where $\mu_{0}$ is a nonvanishing
section of ${\cal L}$ and $\psi (z)\in {\cal O}(U)$.
\\

%%%%%%%%%%%%%%%%%%%%%%%%%%%%%%%%%%%%%%%%%%%%%%%%%%%%%%%%%%%%%%%%%%%

\section{ Generalized Fock-Bargmann representation}

Now we consider the mechanical
system whose quantization is connected with generalization  of
Fock-Bargmann representation.     Let     $(M,\omega)$      be
$2n$-dimensional K\"ahler manifold with fundamental
form $- \omega$ and positively definite K\"ahler metric $g$
to  be  given  in  local complex coordinates  $(z^{\alpha},
z^{\ol\alpha})$ by the formula
\bqq{te15}
-\omega= -\omega^{\alpha\ol\beta}
dz^\alpha \land dz^{\ol\beta}= -i\partial_{\alpha\ol\beta}\Phi\hskip
1pt
dz^\alpha \land dz^{\ol\beta},
\end{equation}
\bqq{te16}
g=g_{\alpha\ol\beta}dz^\alpha dz^{\ol\beta}=
\partial_{\alpha\overline{\beta}}\Phi \hskip 1pt dz^{\alpha}
dz^{\overline{\beta}},
\end{equation}
where $\Phi$ is the K\"ahler potential. As the 2-form
$\omega$ is closed and nondegenerate, it defines the symplectic
structure on   $M$   and   we  can  consider  $(M,\omega)$  as
symplectic manifold and the phase  space  of  some  mechanical
system. The  classical  observables  [{\bf 8}] of such system form Lie
algebra $C^\infty(M)$ with respect to Poisson brackets.

Let as define the K\"ahler polarization $F$ on $TM \otimes
{\bf C}$ in the form
\bqq{te17}
F=\{ X \in TM\otimes {\bf C} \hskip 2pt\vert
X=\xi_{\alpha}V(z^{\alpha}),
\xi_{\alpha}\in C^\infty (M) \},
\end{equation}
where
$V(z^{\alpha})=\omega^{\ol\sigma\alpha}\partial_{\ol\sigma}$
is given by (\ref{te3}). By (\ref{te10}) the function $f$ preserves
the
polarization $F$ if and only if
\bqq{te18}
[V(f),V(z^{\alpha})]=a[f]^{\alpha}_{\mu}V(z^{\mu})
\end{equation}
where according to (\ref{te3}) and (\ref{te15})
$V(f)=\omega^{\mu\ol\nu}
(\partial_{\ol\nu}f\partial_{\mu} - \partial_{\mu}f\partial_
{\ol\nu})$, whence
$$
-   \omega^{\mu\ol\nu}\partial_{\ol\nu} f
\partial_{\mu}\omega^{\alpha\ol\sigma}\partial_{\ol\sigma}
+   \omega^{\mu\ol\nu}\partial_{\mu}f
\partial_{\ol\nu} \omega^{\alpha\ol\sigma}\partial_{\ol
\sigma}
+   \omega^{\alpha\ol\sigma} \partial_{\ol\sigma}
\omega^{\mu\ol\nu} \partial_{\ol\nu} f
\partial_{\mu}
+
$$
$$
\omega^ {\alpha\ol\sigma} \partial_{\ol\nu}\partial_{\ol\sigma}
f \omega^{\mu\ol\nu} \partial_{\mu}
- \omega^{\alpha\ol\sigma}\partial_{\ol\sigma}
\omega^{\mu\ol\nu} \partial_\mu f \partial_{\ol\nu}
- \omega^{\alpha\ol\sigma}\partial_{\ol\sigma}\partial_\mu f
\omega^{\mu\ol\nu} \partial_{\ol\nu} =
- a[f]^\alpha_\mu \omega^{\mu \ol\sigma}\partial_{\ol\sigma}.
$$
Equating the corresponding components of vector  fields  in  left  and
right  parts of the last relation we find
\bqq{te20}
\partial_{\ol\sigma}\omega^{\mu\ol\nu}\partial_{\ol\nu} f +
\omega^{\mu\ol\nu}\partial_{\ol\nu}\partial_{\ol\sigma}f
\partial_\nu=0.
\end{equation}
Differentiating the     equation
$\omega^{\nu\ol\mu}\omega_{\ol\mu\rho} =\delta^\nu_\rho$
we obtain the identity
$\omega^{\nu\ol\mu}\partial_{\ol\gamma}\omega_{\ol\mu\rho} =
-\partial_{\ol\gamma}\omega^{\nu\ol\mu}\omega_{\ol\mu\rho}$.
After this (\ref{te20}) takes the form
\bqq{te21}
\nabla_X \nabla_Y f = 0, \qquad X, Y \in F
\end{equation}
where $\nabla$ denotes   the   covariant  derivation  with  respect
to K\"ahler metric $g$.
By Theorem 1 we find from (\ref{te8})
$\omega=d\alpha$ and from (\ref{te15})
\bqq{te22}
\alpha=-i\hskip 1pt \partial_\alpha\Phi \hskip 1pt dz^\alpha
\end{equation}
modulo to the exact 1-form $d\beta$.

If we choose both $\mu$ and $\lambda$  in  (\ref{te4})
equal  to  nonvanishing
section $\mu_0$, then using (\ref{te6}) we obtain
$$
X<\mu_0,\mu_0>=
i\hbar^{-1}(\alpha(X)-\overline{\alpha(X)})<\mu_0,\mu_0>.
$$
Evaluating this formula on the vector fields
$\partial_\alpha$, $\partial_{\ol\alpha}$, $\alpha=1,...,n$
we find with the help of (\ref{te22})
\bqq{te25}
<\mu_0,\mu_0>= exp(- \hbar^{-1}\Phi)
\end{equation}
up to constant multiplier which we omitted.

Now we determine the sections of Hermitian line bundle ${\cal L}$
which
form the   representation   space   ${\cal   H}$.   Being  covariantly
constant with respect to $D_{X\in F}$ these sections must obey the
equation (\ref{te9}):
$$
D_{V(z^\alpha)} \mu=0, \quad \mu\in\Gamma({\cal L}).
$$
{}From here  we  find with the help of (\ref{te5}),  (\ref{te6}),
(\ref{te17}) and (\ref{te22})
that $\mu =\psi (z)\mu_0$,  where  $\psi  (z)$  is  holomorphic
function on $U \subset {\bf C}^n$.

If $M$ is contractible then using the Hermitian structure
$<,>$ in ${\cal L}$ we  can  define  the
scalar product in ${\cal H}$ by the formula [{\bf 8,12}]
\bqq{te26}
(\mu_1,\mu_2) = \int \psi_1 (z) \overline{ \psi_2 (z)}  <\mu_0,\mu_0>
\omega^n,
\eeq
where $\mu_1=\psi_1(z)\mu_0$, $\mu_2=\psi_2(z)\mu_0$
and $\omega^n$ is n-th
external degree of $\omega$. Using (\ref{te25})  we find
\bqq{te27}
 (\mu_1,\mu_2) = \int \psi_1 (z)\ol { \psi_2 (z)}exp(-\hbar^{-1}\Phi)
\hskip 2pt \omega^n.
\end{equation}
In this case the   representation   Hilbert   space   associated  with
polarization $F$  given  by  (\ref{te17})  can  be  identified  with
Fock space $L_2^{hol}(U,dm)$ of  holomorphic  functions
on  $U\subset{\bf  C}^n$
quadratically integrable     with     the     measure    $dm=$ exp $(-
\hbar^{-1}\Phi)\hskip 1pt \omega^n$.

Let us  consider  the  K\"ahler  space  ${\cal  K}_{2n}$   of
constant
holomorphic curvature $k$ (see for example [{\bf 9}]).
As it is known the space ${\cal K}_{2n}$ is
isometric to the projective space ${\bf CP}^n$ for $k>0$,  to the disk
$D^R_n= \{z\in{\bf C}^n\vert z\ol z<R\}$ for $k<0$ and to ${\bf C}^n$
for $k=0$.  The metric of the  space  ${\cal  K}_{2n}$  in  the  local
complex coordinates is
\bqq{te28}
ds^2= 2g_{\alpha\ol\beta}dz^\alpha dz^{\ol\beta},
\qquad g_{\alpha\ol\beta} =\partial_{\alpha\ol\beta}\Phi=
(A\delta_{\alpha\beta}  -\frac{k}{4}z^{\ol\alpha}z^\beta)A^{-2},
\end{equation}
$$
\Phi=\frac{4}{k}ln\hskip 1pt A,\quad A=1+\frac{k}{4}
\Sigma z^\nu z^{\ol\nu}.
$$

The curve $x(t)$ on K\"ahler manifold $M$ is called $H$-{\it planar}
(or {\it holomorphical planar}) [{\bf 11}] if
it obeys the following equation
$$
\nabla_{\chi} \chi = a(t)\chi +  b(t) J(\chi), \qquad
\chi\equiv \dot x
$$
where $a(t)$,    $b(t)$    are    some   real-valued   functions   and
$J$ is complex structure operator in $TM$.

Let $M$ and $M'$ be two K\"ahler manifolds.
The mapping $f:M\to M'$  is  called  $H$-projective  (see  for
example [{\bf 13}]) if it
transforms $H$-planar
curves of  $M$  into  $H$-planar  curves  of $M'$.

The contravariant components  and  nonvanishing  Christoffel
symbols of the metric (\ref{te28}) in local complex coordinates
are given by the formula
$$
g^{\alpha\ol\beta}= (\delta^{\alpha\beta}
+\frac{k}{4}z^\alpha z^{\ol\beta})A^2 \equiv
-i\omega^{\alpha\ol\beta},\qquad
\Gamma^\alpha_{\beta\gamma}=
-\frac{k}{4}A^{-1}(z^{\ol\alpha}\delta^\beta_\gamma
+z^{\ol\beta}\delta^\alpha_\gamma)=
\overline{\Gamma^{\ol\alpha}_{\ol\beta\ol\gamma}}.
$$
{}From here it is follows that considered metrics are $H$-projectively
flat. Now we find from (\ref{te21})
\bqq{te29}
\partial_{\mu\ol\nu} f +2A^{-1}\partial_{(\ol\mu}A
\partial_{\ol\nu)}f=0.
\end{equation}
After substitution $f= W A^{-1}$ in this equation it takes the form
$$
\partial_{\ol\mu\ol\nu} W=0
$$
whence
$$
W=u_{\ol\alpha}(z)z^{\ol\alpha} + v(z),
$$
where $u_{\ol\alpha}$ and $v$ are arbitrary holomorphic functions.

In [7] the system of observables
\bqq{te32}
\tilde H =\frac{1+z\ol z} {1-z\ol z},\quad N=\frac{z}{1-z\ol z},
\qquad \ol N=\frac{\ol z}{1-z\ol z}
\end{equation}
was considered when quantizing 1-dimensional harmonic  oscilator.
The Ha\-mil\-to\-ni\-an vec\-tor fields $V(\tilde H)$, $V(N)$
and $V(\ol N)$ form the basis of holomorphic isometries Lie algebra
in the  space  ${\cal K}_2$ of holomorphic curvature $k=-4$
(see for example [{\bf 14}]). We use the next system of observables
\bqq{te33}
H =\frac{\Sigma  z^\nu  z^{\ol\nu}}{A},\quad (u_{\ol\alpha} =z^\alpha,
\quad v=0),
\end{equation}
\bqq{te34}
N^\beta =\frac{ z^\beta}{A},\quad (u_{\ol\alpha} =0,
\quad v=z^\beta),
\end{equation}
\bqq{te35}
N^{\ol\beta} =\frac{z^{\ol\beta}}{A},\quad(u_{\ol\alpha}=
\delta_\alpha^\beta,\quad v=0).
\end{equation}
One can easily check that $H$, $N^\alpha$ and
$N^{\ol\alpha}$ are the
solutions of  equation (\ref{te29}).
The Hamiltonian vector fields of this functions define
infinitesimal isometries   in   ${\cal   K}_{2n}$  and  $H$-projective
transformations in  the  flat  K\"ahler manifold ${\bf C}^n$.  Note
that
these isometries do not form Lie algebra.  In  1-dimensional  case
$N^1$, $N^{\ol 1}$  coincide  with  $N,\ol N$  from (24) and $H$ can
be obtained from $\tilde H$ by the linear substitution.  The
using of
$H$ is  more preferable from the point of view of the limit transition
to the flat space ($k=0$). In the limit $k\to 0$
we  obtain  $H=\Sigma\hskip
1pt z^\nu z^{\ol\nu}$,  i.e.  the Hamiltonian of  harmonic  oscillator
(\ref{te12})
written in complex coordinates.

The Hamiltonian  vector  fields  of  the  functions  $H,N^\alpha$  and
$N^{\ol\alpha}$ have the form
$$
T \equiv V(H)=  \omega^{\mu\ol\nu}(\partial_{\ol\nu}H  \partial_\mu
-\partial_\nu H  \partial_{\ol\mu})=  i\hskip 1pt(z^\nu\partial_\nu
-z^{\ol\nu}\partial_{\ol\nu}),
$$
\bqq{te36}
T^\alpha \equiv V(N^\alpha) =-i\hskip 1pt (\frac{k}{4}z^\alpha
z\nu\partial_\nu+\partial_{\ol\alpha}),
\eeq
$$
T^{\ol\alpha} \equiv V(N^{\ol\alpha}) =i\hskip 1pt
(\frac{k}{4}z^{\ol\alpha}
z{\ol\nu}\partial_{\ol\nu}+\partial_\alpha),
$$
Using this  formulae  we  can  calculate the commutators of the vector
fields $T$, $T^\alpha$ and $T^{\ol\alpha}$
\bqq{te37}
[T^\alpha,T^\beta]=0,\qquad
[T^\alpha,T^{\ol\beta}]= i\hskip 1pt \frac{k}{4}(\delta^\alpha_\beta
T+
T^{\alpha\ol\beta}),
\end{equation}
$$
[T^\alpha,T\hskip 1pt]= -i\hskip 1pt T^\alpha,
\qquad [T^{\ol\alpha},T \hskip 1pt]= -iT^{\ol\alpha},
$$
where
$T^{\alpha\ol\beta}=$ $V(z^\alpha    z^{\ol\beta    }A^{-1})=$
$i \: (z^\alpha\partial_\beta- z^{\ol\beta}\partial_{\ol\alpha})$.

The generators $T$, $T^\alpha$ and $T^{\ol\alpha}$ do  not  form  a
Lie
algebra but if we join to them the generator $T^{\alpha\ol\beta}$
then we obtain
\bqq{te38}
[T^\alpha,T^{\beta\ol\gamma}]= i\delta^\alpha_\gamma
T^\beta,\qquad
[T^{\ol\alpha},T^{\beta\ol\gamma}]=-i\delta^{\alpha}_{\beta}
 T^{\ol\gamma},
\end{equation}
$$
[T^{\alpha\ol\beta},T^{\gamma\ol\nu}]=i(\delta^\gamma_\beta
T^{\alpha\ol\nu}-\delta^\nu_\alpha T^{\gamma\ol\beta}).
$$
Because $T^\alpha,$ $T^{\ol\alpha},$ $T^{\alpha\ol\beta}$   are
linearily independent  and $T= \Sigma\hskip 1pt T^{\alpha\ol\alpha}$
from (29), (30) it follows that $T^\alpha,\quad T^{\ol\alpha}$ and
$T^{\alpha\ol\beta}$ form the basis of $n(n+4)$-dimensional (over
${\bf
R}$) Lie  algebra  $ \ell (k)$  which  is  the   Lie   algebra   of
infinitesimal isometries
of  the  space ${\cal K}_{2n}$ preserving the
complex structure.

The Poisson brackets  $\{f,g\}=
\omega^{\alpha\ol\beta}(\partial_\alpha f \partial_{\ol\beta} g-
\partial_{\ol\beta}f \partial_\alpha g)$  of the functions
$H^{\alpha\ol\beta}=z^\alpha z^{\ol\beta}A^{-1}$, $N^\alpha=
z^\alpha A^{-1}$ and $N^{\ol\beta}= z^{\ol\beta}A^{-1}$ are
$$
\{N^\alpha,N^\beta\}=0,\qquad
\{N^\alpha,N^{\ol\beta}\}= i\frac{k}{4}(\delta^\alpha_\beta H+
N^{\alpha\ol\beta})-i\delta^\alpha_\beta,
$$
$$
\{N^\alpha,N^{\beta\ol\gamma}\}= i\hskip 1pt \delta^\alpha_\gamma
N^\beta,
\qquad
\{N^{\ol\alpha},N^{\beta\ol\gamma}\}=  -i\hskip 1pt
\delta^\alpha_\beta N^{\ol\gamma},
$$
$$
\{ N^{\alpha\ol\beta},N^{\gamma\ol\nu}\}=
i(\delta^\alpha_\nu
N^{\gamma\ol\beta} -\delta^\gamma_\beta N^{\alpha\ol\nu}).
$$
Note that if we take the limit $k\to 0$ (\ref{te37}), (\ref{te38})
turns to
$$
[T^\alpha,T^\beta]=0, \qquad
[T^\alpha,T^{\ol\beta}]=0,
$$
\bqq{te40}
[T^\alpha,T^{\beta\ol\gamma}]= i\hskip 1pt \delta^\alpha_\gamma
T^\beta,
\qquad [T^{\ol\alpha},T^{\beta\ol\gamma}]=  -i\hskip 1pt
\delta^\alpha_\beta
T^{\ol\gamma},
\end{equation}
$$
[T^{\alpha\ol\beta},T^{\gamma\ol\nu}]=
i\hskip 1pt (\delta^\alpha_\nu
T^{\gamma\ol\beta} -\delta^\gamma_\beta T^{\alpha\ol\nu}).
$$
and define  the $n(n+4)$-dimensional  Lie  algebra  $\ell (0)$.
The    curve $\ell (k)$   in   the   manifold   of    the
$n(n+4)$-dimensional Lie algebra structures is the deformation of
algebra $\ell (0)$  defined  by the commutation relations
(\ref{te40}) and
containing the  $n$-dimensional  harmonic   oscillator   algebra ${\bf
osc}(n)$ as  the  Lie  subalgebra.  That  is  why  we can consider the
mechanical system with the phase  space  ${\cal  K}_{2n}$,  symplectic
form $\omega$ and the observables $H$, $N^\alpha$ and $N^{\ol\alpha}$
as
the "deformation" of classical $n$-dimensional harmonic oscillator.

Now we  quantize  the  classical  mechanical   systems
obtained in the preceding sections using the polarization $F$ defined
by (\ref{te17}).  We calculate now $a[f]=  \Sigma  a[f]^\nu_\nu$  (see
(\ref{te11}))
for $f=H$, $N^\alpha$ and $N^{\ol\alpha}$. Substituting in
(\ref{te10})
$V(z^\alpha)= \omega^{\ol\nu\alpha}\partial_{\ol\nu}$ instead   of
$X^\alpha$ and $T$, $T^\alpha,T^{\ol\alpha}$ instead of $V(f)$  we
find
using formulae (\ref{te27}) and (\ref{te36})
$$
[T,V(z^\beta)]= iV(z^\alpha),
$$
$$
[T^\alpha,V(z^\beta)]= -i\hskip 1pt
\frac{k}{4}(z^\alpha\delta^\beta_\nu
+z^\beta \delta^\alpha_\nu)V(z^\nu),
$$
$$
[T^{\ol\alpha},V(z^\beta)]=0
$$
whence
$$
a[N^\alpha]=-i \frac{k}{4}z^\alpha (n+1), \qquad a[N]= 0,\qquad
a[H]=in.
$$
After this  from (\ref{te33})-(\ref{te35})
using (\ref{te5}), (\ref{te22}) and (\ref{te26}) we  obtain
the  following  expressions  for  differential
operators in $L^{hol}_2 (U,dm)$ which are  the  quantizations  of  the
observables $H$, $N^\alpha$ and $N^{\ol\alpha}$
\bqq{te45}
{\cal Q}H     \equiv     \hat     H=     \hbar(z^\nu    \partial_\nu
\psi+\frac{n}{2}\psi),
\end{equation}
$$
{\cal Q}N^\alpha  \equiv  \hat  N^\alpha  =-\hbar\frac{k}{4}z^\alpha
z^\nu\partial_\nu \psi+    z^\alpha(1-\frac{\hbar   k}{8}(n+1))
\psi,
$$
$$
{\cal Q}N^{\ol\alpha}        \equiv       \hat       N^{\ol\alpha}
=\hbar\partial_\alpha \psi.
$$

Let $B:{\cal H}\to{\cal H}$ be the  selfadjoint  operator  in  Hilbert
space ${\cal H}$. The set $\sigma (B)= \{\rho \in {\bf R}\vert\hskip
1pt
\exists \hskip 1pt
\mu_\rho \in{\cal  H}:  \hskip 1pt B\mu_\rho=  \rho\mu_\rho  \}$
is  called   the
{\it spectrum} of  the operator $B$.  The number $\rho\in \sigma (B)
\subset
{\bf R}$ and section $\mu_\rho \in  {\cal  H}$
are called the {\it eigenvalue} of $B$ and {\it eigenstate}
with eigenvalue $\rho$.

Let us  consider  the  eigenstate $\psi_E$ of the operator ${\cal Q}H$
with the eigenvalue $E$. Equation (\ref{te45}) yields
$$
\hbar (z^\nu\partial_\nu\psi_E +\frac{n}{2}\psi_E)= E\psi_E
$$
which is equivalent to
$$
z^\nu\partial_\nu\psi_E= (E\hbar^{-1}- \frac{n}{2})\psi_E
$$
whence $\psi_E$ is a homogeneous function of $z$ of degree $l=
(E\hbar^{-1}-
\frac{n}{2})$. Since  $\psi_E$  is holomorphic it follows that $l$
is non-negative integer, so that  the  spectrum of ${\cal Q} H$ is
given by
$$
E_l= (l+ \frac{n}{2})\hbar,\quad l \in \{0\}\cup {\bf N}
$$
and coincides with the spectrum of the $n$-dimensional harmonic
oscillator Hamiltonian (\ref{te12}) (see for example [{\bf 12}]).
\\

%%%%%%%%%%%%%%%%%%%%%%%%%%%%%%%%%%%%%%%%%%%%%%%%%%%%%%%%%%%%%%%%%%%%

\section{ Quantization of K\"ahler manifolds
\newline
admitting $H$-projective mappings}

In this section we  consider  the  quantization  of  K\"ahler
spaces admitting $H$-projective  mappings  onto  another K\"ahler
spaces.  Let
$(M,\omega)$ and  $(M^\prime,   \omega ')$   be    two    K\"ahler
manifolds with fundamental   forms   $-\omega$   and   $-\omega'$.
Let
$\varrho:M \to  M'$ be $H$-projective mapping,  it is well known
that $H$-projective mapping  preserves  the
complex structure. Therefore we can  choose  the local
complex chart $(z^\alpha,z^{\ol\alpha},U)$ in $M$ such that  for  each
point $p  \in  U$  with  coordinates  $(z^\alpha,z^{\ol\alpha})$ its
image $\varrho (p) \in  \varrho  (U)$  has  the  same  coordinates.
The
necessary and  sufficient  condition for the mapping $\varrho:M \to
M'$
to be $H$-projective is expressed with the following  equation [{\bf
13}]
\bqq{te49}
b_{\alpha\ol\beta;\gamma}= 2 \phi'_\alpha g_{\ol\beta\gamma},
\end{equation}
where
$$
b_{\alpha\ol\beta}= e^{2 \phi}g'^{\ol\mu\nu}   g_{\ol\mu\alpha}
g_{\nu\ol\beta}, \qquad b_{\alpha\beta}= b_{\ol\alpha\ol\beta}= 0,
$$
$$
\phi'_{\alpha}= \partial_\alpha \phi' =\partial_\mu \phi
e^{2 \phi}g'^{\mu\ol\nu} g_{\alpha\ol\nu},
$$
$$
J^i_k= J'^i_k,      \qquad      J^\mu_\nu=      -J^{\ol\mu}_{\ol\nu}
=i\delta^\mu_\nu,\qquad J^\mu_{\ol\nu}= J^{\ol\mu}_\nu= 0,
$$
$\phi$ is some function on $U$ and semicolon denotes the
covariant derivation with respect to $g$.

Because of the positive definiteness of $g$ we can define the  complex
frame $\{Z_A,Z_{\ol  A}\}$ which is adapted for the Hermitian
structure
of {\it M} [{\bf 10}]. Then for the frame components of $g$ we have
$$
g_{A \ol B}= \delta^A_B.
$$
Transformations, preserving  this  form  of  $g$,  belong to the
unitary
group $U(n)$ for each point $p \in M$. With the help of such
transformations we can choose the frame $\{Y_A,Y_{\ol B}\}$, so that
\bqq{te51}
g_{A\ol B}=  \delta  ^A_B,
\qquad  b_{A\ol   B}=   \lambda_A \delta^A_B,
\end{equation}
where $\lambda_A={\ol{\lambda}}_A$     are    the    roots    of
the
$\lambda$-matrix $(b- \lambda g)$.  Written  in  the  frame  (see  for
example [{\bf 1,2}]) (\ref{te21}) and (\ref{te49}) have the form
\bqq{te52}
Y_{\ol A} \hskip 2pt Y_{\ol B}f- \sum\limits_S
\gamma_{\ol B S\ol A} Y_S f= 0,
\end{equation}
\bqq{te53}
\delta_{AB}Y_A\lambda_A +\sum\limits_S (\gamma_{\ol S A C}\lambda_S
\delta_{SB}+ \gamma_{S\ol B C}\lambda_S \delta_{SA})= 2\delta_{CB}Y_A
\phi',
\end{equation}
where $\gamma_{\ol S A C}$ $(\gamma_{ S A C}= \gamma_{\ol S \ol A
C})=0$  are Ricci  rotation  coefficients  of  the frame.
\vskip 2pt

Let $\{\theta^A,   \theta^{\ol  A}\},$  $\theta^A,\theta^{\ol  A}\in
T^\ast M$ be the coframe dual to the frame $\{Y^A,Y^{\ol A}\}$.  Then
the connection  form  $\alpha$ in the Hermitian line bundle ${\cal L}$
(see \S {\bf 1}) can be written in the following form
$$
\alpha= -iY_A \Phi \theta^A
$$
as in \S {\bf 2}. From (\ref{te15}) and (\ref{te51}) it follows
$\omega_{A\ol B} =\omega^{A\ol B}= i\delta^A_B$. Then for
$Y_A= \xi^\mu_A \partial_\mu$ we obtain
$$
V(z^\alpha)= -i\hskip 1pt \Sigma\hskip 1pt \xi^\alpha_B\hskip 1pt
 Y_{\ol B}.
$$
If the function $F \in C^\infty (U)$  preserves  polarization  $F$  it
obeys the condition (\ref{te18}) which we can write as the form
$$
[V(f),V(z^\alpha)]= a[f]^\alpha_\mu   V(z^\mu)= \sum\limits_{A,B}
(Y_{\ol A}f   Y_A   \xi^\alpha_B-  Y_A  f  Y_{\ol  A}  \xi^\alpha_B+
\xi^\alpha_A (Y_{\ol A} Y_B f))Y_{\ol B},
$$
and here we  find
$$a[f]= i\zeta^A_\alpha  \sum\limits_B  (Y_{\ol  B}fY_B  \xi^\alpha_A
-Y_B fY_{\ol     B}\xi^\alpha_A+     \xi^\alpha_B    Y_{\ol    B}Y_A
f),                                                            $$
where $\zeta^B_\alpha$   are  components  of  the  inverse  to
$(\xi^\alpha_B)$ matrix:
$\zeta^A_\mu \xi^\nu_A= \delta^\nu_\mu$.

At last we  can evaluate the  differential  operator  ${\cal  Q}f$  in
$L^{hol}_2 (U,dm)$ for function $f\in C^\infty (U)$ obeying
(\ref{te52})
$$
({\cal Q})f \psi\equiv\hat f \psi= \hbar \sum\limits_A Y_{\ol  A}fY_A
\psi- Y_A\Phi   Y_{\ol   A}f  \psi+  f\psi-  \frac{i\hbar}{2}a[f]\psi.
$$
To obtain concrete results we have to  use  specific  expressions  for
Ricci rotation coefficients and frame vector fields.  In particular,
for
4-dimensional K\"ahler manifold admitting $H$-projective mappings
there
are two possibilities
$$
1)\quad\lambda_1= \lambda_2,\qquad 2) \quad \lambda_1 \neq
\lambda_2.
$$
In the first case we have $b= \lambda g$ and,  hence,  $g'=\mu  g$
where $\mu= (\lambda e^{2 \phi})^{-1}$.  Then (\ref{te53})
implies $\mu= const$, and
$H$-projective mappings are only rescalings of the metric.
The second case is more interesting. If $\lambda_1 \neq\lambda_2$ then
from (\ref{te53}) it follows [{\bf 3}]
$$
Y_1 \lambda_1=   Y_{\ol    1}\lambda_2=    Y_2\lambda_1=    Y_{\ol
2}\lambda_1= 0,\qquad \gamma_{1\ol 2 2 }=Y_1 \ln | \lambda_1-
\lambda_2 |,
$$
$$
\gamma_{\ol 1 2 1  }=-Y_2  \ln  |\lambda_1-  \lambda_2|, \qquad
\gamma_{\ol 2 1 1 }=\gamma_{\ol 1 2 2 }.
$$

\bigskip

This work was partially supported by grant 1749
of International Science Foundation and grant RFFI-94-01-01118-a of
Rus\-sian Foun\-da\-tion for Fun\-da\-men\-tal In\-ves\-ti\-ga\-ti\-
ons.

%%%%%%%%%%%%%%%%%%%%%%%%%%%%%%%%%%%%%%%%%%%%%%%%%%%%%%%%%%%%%%
\newpage
\centerline{\bf References}

\medskip

\noindent
[{\bf 1}]  A.V. Aminova: On skew-orthonormal  frame and parallel
symmetric
bilinear form on Riemannian manifolds,
{\it Tensor, N.S.},{\bf 45} (1987), 1-13.

\noindent
[{\bf 2}]  A.V. Aminova: Pseudo-Riemannian manifolds with common
geodesics,
{\it Uspekhi Matematicheskikh Nauk}, {\bf 48} (1993), 107-159.

\noindent
[{\bf 3}]  A.V. Aminova and D.A. Kalinin: $H$-pro\-jec\-ti\-vely
equi\-va\-lent
four-di\-men\-si\-onal Rieman\-nian connections,
{\it IZV. VUZ. Matem.}, No.8 (1994), 11-21.

\noindent
[{\bf 4}]  R. Balbinot, A.El Gradechi, J.-P. Gazeau
and B.Giorgini: Phase space
for quantum elementary systems in anti-de Sitter and Minkowski
spacetimes,
{\it J. Phys. A},{\bf 25} (1992), 1185-1210.

\noindent
[{\bf 5}]  F.A. Berezin: General conception of quantization,
{\it Comm. Math. Phys.},  {\bf 40} (1975), 153-174.

\noindent
[{\bf 6}]  J.-P. Gazeau and V.Hussin: Poicar\'e contraction of
$SU(1,1)$ Fock-Bargmann
structure, {\it J.Phys. A}, {\bf 25} (1992), 1549-1573.

\noindent
[{\bf 7}]  J.-P.  Gazeau and J. Renaud: Lie algorithm for
an interacting $SU(1,1)$
elementary systems and its contraction,
{\it Preprint Universit\'e Paris VII,
PAR-LPTM-92}.

\noindent
[{\bf 8}]  A.A. Kirillov: Geometric quantization. In:"Mo\-dern prob\-
lems of
mathematics" (Itogi Nauki i Tekniki),
{\bf 4}, VINITI SSSR, Moscow, 1985, 141-204 {\sloppy
}

\noindent
[{\bf 9}]  S. Kobayashi and K. Nomizu: Foundation of differential
geometry,
V. II,  Intersci. Publ., N.Y., 1969.

\noindent
[{\bf 10}]  A. Lichnerowicz: Theorie globale des connexiones et des
groupes
d'ho\-lo\-no\-mie,  Cremonese, Roma. 1955.

\noindent
[{\bf 11}]  T. Otsuki and Y. Tashiro: On curves in Kahlerian spaces,
{\it Math. J.  Okayama Univ.}, {\bf  4}(1954), 57-78.

\noindent
[{\bf 12}]  J. Sniatycki: Geometric quantization  and  quantum
mechanics,
Springer, Berlin etc, 1980.

\noindent
[{\bf 13}]  N.S. Sinyukov: Geodesic mappings of Riemannian spaces,
{\it Moscow,  Nauka},  1979.

\noindent
[{\bf 14}]  J.M. Tuynman: Quantization.   Towards  a  comparison
between
methods, {\it J. Math. Phys.}, {\bf 28}(1987). 2829-2840.

%%%%%%%%%%%%%%%%%%%%%%%%%%%%%%%%%%%%%%%%%%%%%%%%%%%%%%%%%%%%%%%%%%%%%-

\end{document}